\documentclass[12pt,leqeq]{article}

\usepackage{amssymb}

\usepackage{amsmath}

\usepackage[all]{xy}

 \newcommand{\Z}{\mathbb{Z}}
 \newcommand{\Agr}{A^\bullet}
 \newcommand{\Bgr}{B^\bullet}
 \newcommand{\Mgr}{M^\bullet}
 \newcommand{\Ngr}{N^\bullet}
 
 \newcommand{\kr}{\mathbf{k}}
 \newcommand{\R}{\mathcal{R}}
 \newcommand{\N}{\mathbb{N}}
 \newcommand{\Real}{\mathbb{R}}
 \newcommand{\Complex}{\mathbb{C}}

\newtheorem{thm}{Theorem}

\newtheorem{lem}[thm]{Lemma}

\newtheorem{prop}[thm]{Proposition} 

\newenvironment{pf}{{\bf Proof.}}{}

\newtheorem{defn}[thm]{Definition}

\newtheorem{exmp}[thm]{Example}

\begin{document}

\setlength{\baselineskip}{16pt}

\title{On $N$-differential graded algebras}

\author{Mauricio Angel\thanks{Work partially
supported by IVIC.}\hspace{.3cm} and Rafael D\'\i az\thanks{Work
partially supported by IVIC.}}

\maketitle

\begin{sloppypar}

\begin{abstract}
We introduce the concept of $N$-differential graded algebras
(N-dga), and study the moduli space of deformations of the
differential of an N-dga. We prove that it is controlled by what
we call the (M,N)-Maurer-Cartan equation.
\end{abstract}

\section*{Introduction}

The goal of this paper is to take the first step towards finding a
generalization of Homological Mirror Symmetry (HMS) \cite{Kon} to
the context of $N$-homological algebra \cite{D-V}. In \cite{Fuk}
Fukaya introduced HMS as the equivalence of the deformation
functor of the differential of a differential graded algebra
associated with the holomorphic structure, with the deformation
functor of an $A_\infty$-algebra associated with the symplectic
structure of a Calabi-Yau variety. This idea motivated us to
define deformation functors of the differential of an
$N$-differential graded algebra. An $N$-dga is a graded
associative algebra $A$, provided with an operator $d:A\to A$ of
degree 1 such that $d(ab)=d(a)b+(-1)^{\bar{a}}ad(b)$ and $d^N=0$.
A nilpotent differential graded algebra (Nil-dga) will be an
$N$-dga for some integer $N\geq 2$. Theorem \ref{Catmon} endows
the category of Nil-differential graded algebras with a symmetric
monoidal structure. We remark that such a monoidal structure
cannot be constructed in a natural way for a fixed $N$ (except for
$N=2$), not even using the $q$-deformed Leibniz rule, see
\cite{Sit}.

\smallskip

In Section 2 we consider deformations of a 2-dga into an $N$-dga.
By deforming 2-dgas one is able to construct a plethora of
examples of $N$-dgas. Roughly speaking Theorem \ref{MC2N} tell us
that a derivation of a 2-dga $d_A+e$ is an $N$-differential iff
\[\left. \begin{array}{cc}
(d_{End}(e)+e^2)^{\frac{N-1}{2}}(d_A+e)=0& \mbox{ for N odd,}\\
(d_{End}(e)+e^2)^{\frac{N}{2}}=0& \mbox{ for N
even.}\end{array}\right.\hspace{1cm}\]

In Section 3 we introduce a general formalism for discrete quantum
mechanics. We introduce these models since they turn out, in a
totally unexpected way, to be relevant in the problem of deforming
an $M$-differential into an $N$-differential with $N\geq M$.
Section 4 contains our main result, Theorem \ref{MCMN} which
provides an explicit identity called the $(M,N)$-Maurer-Cartan
equation that controls deformations of an $M$-complex into an
$N$-complex. The construction of the $(M,N)$-Maurer-Cartan
equation is based on an explicit description of coefficients $c_k$
such that
\[(d_A+e)^N=\sum_{k=0}^{N-1}c_k d_A^k,\]
where $c_k$ depends on $d_A$ and $e$. In Section 5 we define a
functional $cs_{2,2N}$ whose critical points are naturally
determined by the $(2,2N)$-Maurer-Cartan equation.

\smallskip

In conclusion in this paper we introduce the moduli space of
deformations of the differential of an $N$-dga and prove that it
is controlled by a generalized Maurer-Cartan equation. We point
out that our methods and ideas can be applied in a wide variety of
contexts. Examples of $N$-dga's coming from differential geometry
are developed in \cite{AD1}. A $q$-analogue, for $q$ a primitive
$N$-th root of unity, of our main result Theorem \ref{MCMN} is
provided in \cite{AD2}. In \cite{AD3} we state an $N$-generalized
Deligne's principle and use the constructions of this paper to
study $A_\infty$-algebras of depth $N$.

\section{$N$-differential graded algebras and modules}

Throughout this paper we shall work with the abelian category of
$\kr$-modules over a commutative ring $\kr$ with unit. We will
denote by $\Agr$ $\Z$-graded $\kr$-modules $\oplus_{i\in\Z}A^i$.
We let $\bar{a}\in\Z$ denote the degree of the element $a\in
A^{\bar{a}}$. The following definition is taken from \cite{Kap}.

\begin{defn}
Let $N\geq 1$ an integer. An {\bf $N$-complex} is a pair $(\Agr,
d)$, where $\Agr$ is a $\Z$-graded object and  $d:\Agr\to \Agr$ is
a morphism of  degree 1 such that $d^N=0$.
\end{defn}

Clearly an $N$-complex is a $P$-complex for all $P\geq N$. If
$\kr$ is a field, then an $N$-complex $(\Agr,d)$ is referred as an
{\bf $N$-differential graded vector space} ($N$-dgvect). An
$N$-complex $(\Agr,d)$ such that $d^{N-1}\neq 0$ is said to be a
{\bf proper $N$-complex}. Let $(\Agr,d_A)$ be an $M$-complex and
$(\Bgr,d_B)$ be an $N$-complex, a {\bf morphism} $f:(\Agr,d_A)\to
(\Bgr,d_B)$ is a morphism $f:\Agr\to\Bgr$ of $\kr$-modules such
that $d_Bf=fd_A$.

\begin{lem}\label{inso}
Let $(\Agr,d_A)$ be a proper $M$-complex, $(\Bgr,d_B)$ be a proper
$N$-complex and $f:(\Agr,d_A)\to (\Bgr,d_B)$ be a morphism, then
(1) If $Ker(f)=0$, then $M\leq N$; (2) If $Im(f)=\Bgr$, then
$M\geq N$ and (3) If $Ker(f)=0$ and $Im(f)=\Bgr$, then $M=N$.
\end{lem}
\begin{pf}
1. Assume that $N<M$ and let $a\in\Agr$ then
$f(d_A^N(a))=d_B^N(f(a))=0$. This implies that $d_A^N(a)\in
Ker(f)=0$, and therefore $d_A^N(a)=0$ which is in contradiction
with the fact that $(\Agr,d_A)$ is a proper $M$-complex. The proof
of 2. is analogous to 1., 3. follows from 1. and 2.$\blacklozenge$
\end{pf}

\begin{exmp}\label{ej1}
Consider $V=\mathbb{C}\!<\!e_1,e_2,e_3\!>$ the complex vector
space generated by $e_1,e_2,e_3$. We endow $V$ with a
$\Z$-graduation declaring $\bar{e_1}=0$, $\bar{e_2}=1$ and
$\bar{e_3}=2$. Define the linear map $d:V\to V$ on generators by
\[d(e_1)=e_2,\hspace{.5cm}d(e_2)=e_3,\hspace{.5cm}\text{and}\hspace{.2cm}d(e_3)=0.\]
$(V,d)$ is a proper 3-complex.
\end{exmp}

\begin{defn}
Let $(\Agr,d)$ be an $N$-complex, we say that an element $a\in
A^i$ is {\bf p-closed} if $d^p(a)=0$ and is {\bf p-exact} if there
exists an element $b\in A^{i-N+p}$ such that $d^{N-p}(b)=a$, for
$1\leq p<N$ fixed. The {\bf cohomology groups} of are the
$\kr$-modules
\[_pH^{i}(A)=\frac{Ker\{d^{p}:A^{i}\to A^{i+p}\}}{Im\{d^{N-p}:A^{i-N+p}\to A^{i}\}},\]
where $i\in\Z,\ p=1,2,...,N\!-\!1$. We set $_kH^*(A)=0$ for $k\geq
N$.
\end{defn}

Notice that a 2-complex $\Agr$ is just a complex in the usual
sense and in this case $p$ is necessarily equal to 1 and
$_1H^{i}(A)$ agrees with $H^{i}(A)$ for all $i\in\Z$.

\begin{defn}\label{Ndga}
\begin{enumerate}
\item[(a)]Let $N\geq 1$ an integer. An {\bf $N$-differential graded algebra} or $N$-dga over $\kr$, is
a triple $(\Agr,m,d)$ where $m:A^k\otimes A^l\to A^{k+l}$ and
$d:A^k\to A^{k+1}$ are $\kr$-modules homomorphisms satisfying
\begin{enumerate}
\item[1)] The pair $(\Agr,m)$ is a graded associative algebra.

\item[2)] For all $a,\ b\in\Agr$, $d$ satisfies
the graded Leibniz rule $d(ab)=d(a)b+(-1)^{\bar{a}}ad(b)$.

\item[3)] $d^{N}=0$, i.e., $(\Agr,d)$ is an $N$-complex.
\end{enumerate}

\item[(b)]A {\bf nilpotent differential graded algebra} (Nil-dga) is an $N$-dga
for some integer $N\geq 2$.
\end{enumerate}
\end{defn}

A 1-dga is a graded associative algebra. A 2-dga is a differential
graded algebra.

\begin{lem}\label{nho} Let $(\Agr,m,d)$ be an $N$-dga, then
if $a$ is p-closed and $b$ is q-closed then $ab$ is
(p+q-1)-closed.
\end{lem}
\begin{pf}
The Lemma follows from the identity
\[d^{n}(ab)=\sum_{i=0}^{n}\left\{\begin{array}{c}n\\i\end{array}\right\}_{\bar{a}}d^{i}(a)d^{n-i}(b),\]
where
$\left\{\begin{array}{c}n\\0\end{array}\right\}_{\bar{a}}=(-1)^{\bar{a}}$,
and for $j\geq 1$,
\(\left\{\begin{array}{c}n+1\\j\end{array}\right\}_{\bar{a}}=
\left\{\begin{array}{c}n\\j-1\end{array}\right\}_{\bar{a}}+(-1)^{\bar{a}+j}
\left\{\begin{array}{c}n\\j\end{array}\right\}_{\bar{a}}.\)

When $n=p+q-1$, since $d^i(a)=0$ for $i\geq p$, we only consider
the case $i<p$, then $n-i=p+q-1-i>q-1$ and $d^{n-i}(b)=0$, because
$d^j(b)=0$ for $j\geq q$. Thus so either $d^i(a)=0$ or
$d^{n-i}(b)=0$ for all $i$, and we have $ab$ is $(p+q-1)$-closed.
$\blacklozenge$
\end{pf}

\begin{defn}
Let $(\Agr,m_A,d_A)$ be an $M$-dga and $(\Bgr,m_B,d_B)$ be an
$N$-dga. A {\bf morphism} $f:\Agr\to\Bgr$ is a linear map such
that $fm_A=m_B(f\otimes Id)+m_B(Id\otimes f)$ and $d_Bf=fd_A$.
\end{defn}

A morphism $f:\Agr\to\Bgr$ such that $f(A^i)\subset B^{i+k}$ is
said to be a morphism of degree $k$. A pair of morphisms
$f,g:\Agr\to\Bgr$ of $N$-dga are homotopic, if there exist
$h:\Agr\to\Bgr$ of degree $N-1$ such that
\[f-g=\sum_{i=0}^{N-1}d_B^{N-1-i}hd_A^i.\]
We remark that if two morphisms $f,g:\Agr\to\Bgr$ of Nil-dga are
homotopic then they induce the same maps in cohomology.

\smallskip

Let $(\Agr,m_A,d_A)$ and $(\Bgr,m_B,d_B)$ be an $M$-dga and an
$N$-dga, respectively. Define $d_{\scriptsize{A\otimes
B}}=d_A\otimes Id+Id\otimes d_B$, the identity
\[d_{\scriptsize{A\otimes B}}^n(a\otimes b)=
\sum_{k=0}^n(-1)^{\bar{a}(n-k)}d_A^k(a)\otimes
d_B^{n-k}(b)\hspace{.3cm}\mbox{ implies, }\]
\begin{prop}\label{most}
The triple $(\Agr\otimes\Bgr,m_{\scriptsize{A\otimes
B}},d_{\scriptsize{A\otimes B}})$ is an $(M\!+\!N\!-\!1)$-dga,
where $m_{\scriptsize{A\otimes B}}=m_A\otimes m_B$.
\end{prop}

\begin{exmp}
Let $(V,d)$ be the 3-complex of in Example \ref{ej1}. On the space
$V\otimes V^*$ consider the base given by $E_{ij}=e_i\otimes
e_j^*$, $i,j=1,2,3$, and define
\[D(E_{ij})=E_{(i+1)j}+(-1)^{i+j}E_{i(j-1)},\]
by Proposition \ref{most} and since $D^4(E_{13})\neq 0$, then
$(V\otimes V^*,D)$ is a proper 5-dga.
\end{exmp}

\begin{thm}\label{Catmon}
The category Nil-dgvect is a symmetric monoidal category. Nil-dga
is the category of monoids in Nil-dgvect. Nil-dga inherits a
symmetric monoidal structure from Nil-dgvetc.
\end{thm}

Let $V^\bullet$ be an $N$-dga. By Proposition \ref{most},
$(V^\bullet)^{\otimes 2}$ is a $(2N\!-\!1)$-dga,
$(V^\bullet)^{\otimes 3}$ is a $(3N\!-\!2)$-dga and in general
$(V^\bullet)^{\otimes k}$ is a $[k(N\!-\!1)\!+\!1]$-dga.

\begin{defn}
Let $(\Agr,m_A,d_A)$ be an $N$-dga and $\Mgr$ a graded
$\kr$-module. Let $K\geq 2$ an integer. A {\bf $K$-differential
graded module} ($K$-dgm) over $(\Agr,m_A,d_A)$, is a triple
$(\Mgr,m_M,d_M)$ with $m_M:A^k\otimes M^l\to M^{k+l}$ and
$d_M:M^k\to M^{k+1}$, $\kr$-modules morphisms satisfying the
following properties
\begin{enumerate}
\item For all
$a,b\in\Agr$ and $m\in\Mgr$, $m_M(a,m_M(b,m))=m_M(m_A(a,b),m)$, .
If no confusion arises, we denote $m_M(a,m)$ by $am$.

\item For all
$a\in\Agr$ and $m\in\Mgr$,
$d_M(am)=d_A(a)m+(-1)^{\bar{a}}ad_M(m)$.

\item The pair $(\Mgr,d_M)$ is a $K$-complex, $d_M^{K}=0$.
\end{enumerate}
\end{defn}

Let $(\Mgr,m_M,d_M)$ be a $K$-dgm and $(\Ngr,m_N,d_N)$ be an
$L$-dgm both over an $N$-dga $(\Agr,m_A,d_A)$. A morphism
$f:\Mgr\to\Ngr$ of degree $k$ is a linear map such that
$f(m_M(a,b))=(-1)^{\bar{a}\bar{f}}m_N(a,f(b))$ and
$d_M(f(b))=f(d_N(b))$, for all $a\in\Agr$ and $b\in\Mgr$. Now let
$(\Mgr,m_M,d_M)$ be a $K$-dgm over an $M$-dga $(\Agr,m_A,d_A)$ and
$(\Ngr,m_N,d_N)$ an $L$-dgm over an $N$-dga $(\Bgr,m_B,d_B)$. The
triple $(M\otimes N,m_{\scriptsize{N\otimes
M}},d_{\scriptsize{M\otimes N}})$ turns out to be a
$(K\!+\!L\!-\!1)$-dgm over $(A\otimes B,m_{\scriptsize{A\otimes
B}},d_{\scriptsize{A\otimes B}})$, where $m_{\scriptsize{M\otimes
N}}$ and $d_{\scriptsize{M\otimes N}}$ are defined as before.

\begin{defn}
The space of endomorphisms of degree $k$ of $\Mgr$ is
\(End^{k}(M)=\prod_{i\in\Z}Hom(M^i,M^{i+k})\), this is,
$End^{k}(M)$ consist of maps $f:\Mgr\to\Mgr$ of degree $k$ which
are linear in regard to the action of $\Agr$ but which does not
satisfy necessarily the relation $d_Mf=(-1)^{\bar{f}}fd_M$.
\end{defn}

There are operators $\circ_M:End(M)\otimes \Mgr\to \Mgr$ and
$\circ_E: End(M)\otimes End(M)\to End(M)$. Similarly to
Proposition \ref{most}, Proposition \ref{Nend} below provides the
natural algebraic structure on $End(M)$,

\begin{prop}\label{Nend}
Define $d_{End}(f)\!\!:=\!\!d_M(f)\!-\!(-1)^{\bar{f}}\!f(d_M)$,
for $f\!\!\in\!\! End(M)$. The triple $(End(M),\circ_E,d_{End})$
is a $(2N\!-\!1)$-dga, and $(\Mgr,\circ_M,d_M)$ is an $N$-dgm over
$(End(M),\circ_E,d_{End})$.
\end{prop}
\begin{pf}
Associativity of $\circ_E$ follows from the associativity
morphisms composition. The Leibniz rule for $d_{End}$ is a
consequence of the Leibniz rule for $d_M$. From the definition of
$d_{End}$ we obtain the identity
\[d_{End}^n(f)=
\sum_{k=0}^n(-1)^{\bar{f}(n-k)} d_M^k\circ f\circ d_M^{n-k}\]
which can be proved by induction and holds for all $n\geq 1$. Let
$n=2N-1$ if $k<N$ then $N-1<n-k$ and thus $d_M^{n-k}=0$. Similarly
if $n-k<N$ then $d_M^k=0$.$\blacklozenge$
\end{pf}

\section{Deformation theory of 2-dgas into N-dgas}

Let $\kr$ a field and consider the category {\bf Artin} of finite
dimensional local $\kr$-algebras. If $\R\in Ob(\mathbf{Artin})$
with maximal ideal $\R_+$ then $\kr\cong \R/\R_+$ ($\R=\kr[[t]]$
and $\R_+=t\kr[[t]]$ are examples to keep in mind). Since
$\kr\cong \R/\R_+$ then $\R\cong \kr\oplus\R_+$ as vector spaces.
We study deformation theory using the formalism which considers
deformations as functors from Artin algebras to Sets for later
convenience.

\begin{defn}\label{def}
Let $\Agr$ be an $M$-dga, an {\bf $N$-deformation} of $\Agr$ over
$\R$ is an $N$-dga $\Agr_\R$ over $\R$, with $N\geq M$, such that
$\Agr_\R/\R_+\Agr_\R$ is isomorphic to $\Agr$ as $N$-dga. Two
$N$-deformations $\Agr_\R$ and $\Bgr_\R$ are said to be {\bf
isomorphic} if there exists an isomorphism
$\Phi:\Agr_\R\to\Bgr_\R$ of $N$-dgas such that the induced
isomorphism $\bar{\Phi}:\Agr_\R/\R_+\Agr_\R\to\Bgr_\R/\R_+\Bgr_\R$
satisfies $i_B\bar{\Phi}=i_A$, where $i_A$ and $i_B$ are the
isomorphism $i_A:\Agr_\R/\R_+\Agr_\R\to\Agr$ and
$i_B:\Bgr_\R/\R_+\Bgr_\R\to\Agr$.
\end{defn}

The core of Definition \ref{def} is to require that $d_{A_\R}$
reduces to $d_A$, and $m_{A_\R}$ reduces to $m_A$ under the
natural projection $\pi: \Agr_\R\to \Agr_\R/\R_+\Agr_\R\cong
\Agr$. Assume that $\Agr_\R=\Agr\otimes\R$ as graded algebras. We
have the following decomposition
\[\Agr_\R=\Agr\otimes\R \\
=\Agr\otimes(\kr\oplus\R_+) \\
=(\Agr\otimes\kr)\oplus(\Agr\otimes\R_+) \\
=\Agr\oplus(\Agr\otimes\R_+). \] Thus, since $d_{A_\R}$ reduces to
$d_A$ under the projection $\pi$, we must have
\[\fbox{$d_{A_\R}=d_A+e$}\]
where $e\in Der(\Agr\otimes\R_+)$ has degree 1. Moreover, the fact
that $d_{A_\R}^N=0$ implies that $e$ is required to satisfy an
identity which we call the $(M,N)$-Maurer-Cartan equation. Next
proposition is well known and considers the classical case, that
is, the $(2,2)$-Maurer-Cartan equation.

\begin{prop}\label{MC2}
Let $\Agr$ be a 2-dga and $\Agr_\R=\Agr\otimes\R$ be a
2-deformation over $\R$, $d_{A_\R}=d_A+e$ where $e\in
Der(\Agr\otimes\R_+)$, then $e$ satisfies the (2,2)-Maurer-Cartan
equation given by
\[d_{End}(e)+e^2=0.\]
\end{prop}
\vspace{-.8cm}
\begin{eqnarray*}
\mbox{ {\bf PROOF.} We have } \hspace{.8cm}d_{A_\R}^{2}(a)&=&(d_A+e)(d_A+e)(a)\hspace{6.3cm}\\
&=&d_A^2(a)+d_A(e(a))+e(d_A(a))+e^2(a)\\
&=&d_{End}(e)(a)+e^2(a),\hspace{.3cm}\mbox{  for all $a\in\Agr.\blacklozenge$ }\\
\end{eqnarray*}

\vspace{-.5cm}

Suppose that $N=2k+n$, $n\in\{0,1\}$ and $k\in\N$, then
\[d_{A_\R}^N=d_{A_\R}^{2k+n}=(d_{A_\R}^2)^kd_{A_\R}^n=(d_{End}(e)+e^2)^kd_{A_\R}^n,\ \mbox{ thus }\]

\begin{thm}\label{MC2N}
Let $\Agr$ be a 2-dga. $\Agr_\R=\Agr\otimes\R$ is an
$N$-deformation over $\R$ with $d_{A_\R}=d_A+e$ where $e\in
Der(\Agr\otimes\R_+)$ of degree 1, iff $e$ satisfies
\[\left.\begin{array}{cc}
(d_{End}(e)+e^2)^{\frac{N-1}{2}}(d_A+e)=0& \mbox{ for N odd,}\\
(d_{End}(e)+e^2)^{\frac{N}{2}}=0& \mbox{ for N
even.}\end{array}\right.\]
\end{thm}

Theorem \ref{MC2N} can be easily extended to study deformations of
the differential of a 2-dgm $\Mgr$ over a 2-dga $\Agr$ as follows

\begin{thm}\label{MC2N-Mod}
Let $\Mgr$ be a 2-dgm over a 2-dga $\Agr$. Then
$\Mgr_\R=\Mgr\otimes\R$ is an $N$-deformation over $\R$ with
$d_{A_\R}=d_A+e$ where $e\in End(\Mgr\otimes\R_+)$ has degree 1,
iff $e$ satisfies
\[\left.\begin{array}{cc}
(d_{End}(e)+e^2)^{\frac{N-1}{2}}(d_M+e)=0& \mbox{ for N odd,}\\
(d_{End}(e)+e^2)^{\frac{N}{2}}=0& \mbox{ for N
even.}\end{array}\right.\]
\end{thm}

Let $M$ be a 3-dimensional smooth manifold. The space
$(\Omega^{\bullet}(M),d)$ of differential forms on $M$ is a
differential graded algebra with $d$ the de Rham differential. Let
$\pi:E\to M$ be a vector bundle, the space
$(\Omega^{\bullet}(M,E),d_E)$ of $E$-valued forms is a
differential graded module over $(\Omega^{\bullet}(M),d)$, where
$d_E$ is the differential induced by $d$. Let $A\in\Omega^{1}(M)$
and consider the endomorphism $e_A$ induced by $A$, defined by
$e_A(\omega)=A\wedge\omega$ for all
$\omega\in\Omega^{\bullet}(M,E)$. The pair
$(\Omega^{\bullet}(M,E),d_E+e_A)$ is a 4-dgm for any $A$.
Moreover, according to Theorem \ref{MC2N-Mod}
$(\Omega^{\bullet}(M,E),d+e_A)$ is a 3-dgm if and only if for all
$\omega$
\[d_{End}(e_A)(d+e_A)\omega=0.\]
Since $d_{End}(e_A)(d+e_A)$ is an operator of degree 3, the
identity $d_{End}(e_A)(d+e_A)\omega=0$ holds for any $k$-form
$\omega$, $k\geq 1$. Thus $(\Omega^{\bullet}(M,E),d+e_A)$ is a
3-dgm if and only if for any 0-form $\omega$
\[d_{End}(e_A)(d+e_A)\omega=d(A)\wedge(d_E(\omega)+A\wedge\omega)=0.\]

Similarly, it is easy to deduce from Theorem \ref{MC2N-Mod} that
if $M$ is an $n$-dimensional smooth manifold and $n<m$, then
$(\Omega^{\bullet}(M,E),d+e_A)$ is a $m$-complex. Let now $M$ be a
$2n$-dimensional smooth manifold. Using local coordinates the
2-form $d_{End}(e_A)$ can be written as $F_{ij}dx^i\wedge dx^j$
where $F_{ij}=\partial_i A_j-\partial_j A_i$. Furthermore,
\[(F_{ij}dx^i\wedge dx^j)^n=\biggl(\sum_{\alpha\in P(2n)}\prod_{i=1}^n
sign(\alpha)F_{a_i,b_i}\biggr)dx^1\wedge...\wedge dx^{2n},\] where
$P(2n)$ is the set of ordered pairings of $[2n]=\{1,...,2n\}$.
Recall that a ordered pairing $\alpha\in P(2n)$ is a sequence
$\{(a_i,b_i)\}_{i=1}^n$ such that $[2n]=\bigsqcup_{i=1}^n
\{a_i,b_i\}$ and $a_i<b_i$. By Theorem \ref{MC2N-Mod},
$(\Omega^{\bullet}(M,E),d+e_A)$ is a $2n$-complex if and only if
the 2-form $F_{ij}dx^i\wedge dx^j$ satisfies
\[\sum_{\alpha\in P(2n)}sign(\alpha)\prod_{i=1}^n F_{a_i,b_i}=0.\]

\smallskip

Let $M$ be a complex manifold and consider the differential graded
algebra $(\Omega(M),\wedge,\bar{\partial})$, where
$\bar{\partial}$ is the Dolbault differential. Let $\pi:E\to M$ be
a complex vector bundle, we consider $\Omega(M,E)$ the forms with
values in $E$. Recall that a holomorphic structure on $E$ is given
by a left differential graded module structure
$(\Omega(M,E),\wedge_E,\bar{\partial_E})$ over the 2-dga
$(\Omega(M),\wedge,\bar{\partial})$. Suppose that on
$(\Omega(M,E),\wedge_E,\bar{\partial_E})$ there is a left
$N$-differential graded module structure over the 2-dga
$(\Omega(M),\wedge,\bar{\partial})$, then in this case we say that
$E$ carries an {\bf $N$-holomorphic structure}.

\section{Discrete quantum theory}

Generally speaking the following data constitute the basic set up
for a (non-relativistic) {\bf quantum mechanical system}: A finite
dimensional Riemannian manifold $M$ which is thought as the
configuration space of the quantum system; A Lagrangian function
$L:TM\to\Real$ which assigns weights to points in phase space.

\smallskip

Associated to this data is the Hilbert space $\mathcal{H}$ of
quantum states which is usually taken to be $L^2(M)$, the space of
square integrable functions on $M$. The dynamics of the quantum
system is determined by operators $U_t:\mathcal{H}\to\mathcal{H}$,
where $t\in\Real$ represents time. The kernel $\omega_t$ of $U_t$
is such that
\[(U_tf)(y)=\int_M \omega_t(y,x)f(x)dx.\]
The key insight of Feynman is that $\omega_t(y,x)$ admits an
integral representation
\[\omega_t(y,x)=\int e^{i\int_0^t L(\gamma,\dot{\gamma})dt}D(\gamma).\]
The integral above runs over all paths $\gamma:[0,t]\to M$ such
that $\gamma(0)=x$ and $\gamma(t)=y$. Making rigorous sense of
this integral is the main obstacle in turning quantum mechanics a
fully rigorous mathematical theory. Recall that a directed graph
$\Gamma$ is given by: $i)$ A set $V_\Gamma$ called the set of
vertices, $ii)$ A set $E_\Gamma$ called the set of edges and
$iii)$ A map $(s,t):E_\Gamma\to V_\Gamma\times V_\Gamma$.
Following the pattern above, one may define a {\bf discrete
quantum mechanical system} as being given by the following data
\begin{enumerate}
\item A directed graph $\Gamma$ (finite or infinite) which plays the role
of configuration space.

\item A map $L:E_\Gamma\to\Real$ called the
Lagrangian map of the system.
\end{enumerate}
The associated Hilbert space is $\mathcal{H}=\Complex^{V_\Gamma}$.
The operators $U_n:\mathcal{H}\to\mathcal{H}$, where $n\in\Z$
represents discretized time are given by
\[(U_nf)(y)=\sum_{x\in V_\Gamma}\omega_n(y,x)f(x),\]
where the discretized kernel $\omega_n(y,x)$ admits the following
representation
\[\omega_n(y,x)=\sum_{\gamma\in P_n(\Gamma,x,y)}\prod_{e\in\gamma}e^{iL(e)}.\]
Here $P_n(\Gamma,x,y)$ denotes the set of length $n$ paths in
$\Gamma$ from $x$ to $y$, i.e., sequences $(e_1,\cdots,e_n)$ of
edges in $\Gamma$ such that $s(e_1)=x$, $t(e_i)=s(e_{i+1}),\
i=1,...,n-1$ and $t(e_n)=y$.

\smallskip

In Section 4 we show that the generalized Maurer-Cartan equation
controlling deformations of $N$-dgas is determined by the kernel
of a discrete quantum mechanical system $L$ which we proceed to
introduce. Let us first explain our notation and conventions which
generalize those introduced in \cite{DiPa}.

\smallskip

For $s=(s_1,...,s_n)\in \N^n$ we set $l(s)=n$, the length of the
vector $s$, and $|s|=\sum_i{s_i}$. For $1\leq i <n,\ s_{>i}$
denotes the vector given by $s_{>i}=(s_{i+1},...,s_n)$, for
$1<i\leq n,\ s_{<i}$ stands for $s_{<i}=(s_1,...,s_{i-1})$, we
also set $s_{>n}=s_{<1}=\emptyset$. $\N^{(\infty)}$ denotes the
set $\bigsqcup_{n=0}^{(\infty)}\N^n$, where by convention
$\N^{(0)}=\{\emptyset\}$.

\smallskip

We define maps $\delta_i,\ \eta_i:\N^n\to \{0,1\}$, for $1\leq
i\leq n$, as follows
\[\delta_i(s)=\left\{\begin{array}{cc}
1 & \mbox{ if $s_i=0$, }\\
0 & \mbox{otherwise. }\end{array}\right. \hspace{2cm}
\eta_i(s)=\left\{\begin{array}{cc}
1 & \mbox{ if $s_i\geq 1$, }\\
0 & \mbox{otherwise. }\end{array}\right.\]

\smallskip

For an $M$-dga $\Agr$ and $e\in End(\Agr)$ and $s\in\N^n$ we
define $e^{(s)}=e^{(s_1)}\cdots e^{(s_n)}$, where
$e^{(l)}=d_{End}^{l}(e)$ if $l\geq 1$, $e^{(0)}=e$ and
$e^{\emptyset}=1$. In the case that $e_a\in End(\Agr)$ is given by
\[e_a(\phi)=a\phi,\hspace{.3cm}\text{for $a\in A^1$ fixed and all $\phi\in\Agr$,}\]
then $e_a^{(l)}=d_{End}^{l}(e_a)$ reduces to
$e_a^{(l)}=e_{d^{l}(a)}$, thus
\[e_a^{(s)}=e_a^{(s_1)}\cdots e_a^{(s_n)}=e_{d^{s_1}(a)}\cdots e_{d^{s_n}(a)}.\]

Where $[k]$ denotes the set $\{1,2,...,k\}$. For $N\in\N$ we
define $E_N=\{s\in\N^{(\infty)}:|s|+l(s)\leq i\}$ and for $s\in
E_N$ we define $N(s)\in\Z$ by $N(s)=N-|s|-l(s)$.

\smallskip

We introduce the discrete quantum mechanical system $L$ by
\begin{enumerate}
\item $V_{L}=\N^{(\infty)}$.
\item There is a unique directed edge in $L$ from vertex $s$ to
$t$ if and only if $t\in\{(0,s),s,(s+e_i)\}$ where
$e_i=(0,..,\underset{\scriptsize{i-th}}{\underbrace{1}},..,0)\in\N^{l(s)}$,in
this case we set $source(e)=s$ and $target(e)=t$.

\item Edges in $L$ are weighted according to the following table
\smallskip
\begin{center}
\begin{tabular}{|l|l|l|}\hline $source(e)$    &      $target(e)$ &
$weight(e)$      \\  \hline
$s$         &      $(0,s)$      & $1$                 \\
$s$         &      $s$          & $(-1)^{|s|+l(s)}$    \\
$s$         &     $(s+e_i)$    & $(-1)^{|s_{<i}|+i-1}$ \\
\hline
\end{tabular}
\end{center}
\end{enumerate}

The set $P_N(\emptyset,s)$ consists of all paths
$\gamma=(e_1,...,e_N)$, such that $source(e_1)=\emptyset$,
$target(e_N)=s$ and $source(e_{l+1})=target(e_l)$. For $\gamma\in
P_N(\emptyset,s)$ we define the weight $\omega(\gamma)$ of
$\gamma$ as
\[\omega(\gamma)=\prod_{i=1}^{N}\omega(e_i).\]

\section{The (M,N)-Maurer-Cartan equation}

\begin{lem}\label{Lem23}
Let $\Agr$ be an $M$-dga and $\R\in Ob({\bf Artin)}$. We define
$d_{A_\R}=d_A+e$ where $e\in Der(\Agr\otimes\R_+)$ has degree 1,
then
\[(d_{A_\R})^N=\sum_{s\in E_N}c(s,N)e^{(s)}d_A^{N(s)},\]
where the coefficient $c(s,N+1)$ is equal to
\begin{equation}
\delta_{1}(s)c(s_{>1},N)+(-1)^{|s|+l(s)}c(s,N)+
\sum_{i=1}^{l(s)}\eta_{i}(s)(-1)^{|s_{<i}|+i-1}c(s-e_i,N),\label{eq0re}\end{equation}
and $c(\emptyset,1)=c(0,1)=1$.
\end{lem}
\begin{pf}
We use an induction on $N$. For $N=1$, since $E_1=\{s=\emptyset,\
s=0\}$
\begin{eqnarray*} d_{A_\R}&=&\sum_{s\in
E_1}c(s,1)e^{(s)}d_A^{1(s)}=
c(\emptyset,1)e^{(\emptyset)}d_A^{1-|\emptyset|-l(\emptyset)}+c(0,1)e^{(0)}d_A^{1-|0|-l(0)}\\
&=&c(\emptyset,1)d_A+c(0,1)e.\\
\end{eqnarray*}

\vspace{-.8cm}

Suppose our formula holds for $N$ and let us check it for $N+1$
\begin{eqnarray}
(d_{A_\R})^{N+1}&=&(d_A+e)(d_{A_\R})^N \nonumber\\
&=&(d_A+e)\biggl(\sum_{s\in E_N}c(s,N)e^{(s)}d_A^{N(s)}\biggr)\nonumber\\
&=&d_A\biggl(\sum_{s\in E_N}c(s,N)e^{(s)}d_A^{N(s)}\biggr)+
e\biggl(\sum_{s\in E_N}c(s,N)e^{(s)}d_A^{N(s)}\biggr)\nonumber\\
&=&\sum_{s\in E_N}c(s,N)d_A(e^{(s)}d_A^{N(s)})+\sum_{s\in
E_N}c(s,N)ee^{(s)}d_A^{N(s)}. \label{eq1}
\end{eqnarray}
Consider the second term of the right hand side of (\ref{eq1})
\begin{eqnarray}
\sum_{s\in E_N}c(s,N)ee^{(s)}d_A^{N(s)}&=&
\sum_{s\in E_N}c(s,N)e^{(0)}e^{(s)}d_A^{N(s)} \nonumber\\
&=&\sum_{\begin{subarray}{c} t\in E_{N+1} \\
t_1=0\end{subarray}}
c(t_{>1},N)e^{(t)}d_A^{N-|t_{>1}|-l(t_{>1})} \label{eq2}\\
&=&\sum_{s\in E_{N+1}}\delta_{1}(s)c(s_{>1},N)e^{(s)}
d_A^{N(s)+1}. \label{eq3}
\end{eqnarray}
In (\ref{eq2}) we put $t=(0,s)$ thus $|t|=|s|$ and $l(t)=l(s)+1$
and (\ref{eq3}) is obtained by rewriting and changing $t$ by $s$.

Now consider the first term of the right hand side of (\ref{eq1})
\begin{eqnarray}
\sum_{s\in
E_N}c(s,N)d_A(e^{(s)}d_A^{N(s)})&=&\sum_{\begin{subarray}{l} s\in
E_N\\1\leq i \leq
l(s)\end{subarray}}(-1)^{|s_{<i}|+i-1}c(s,N)e^{(s+e_i)}d_A^{N(s)} \nonumber \\
& &\hspace{-1cm}+\sum_{s\in E_N}(-1)^{|s|+l(s)}c(s,N)e^{(s)}d_A^{N(s)+1} \label{eq4}\\
& &\hspace{-3.5cm}=\sum_{t\in E_{N+1}}\sum_{\begin{subarray}{c}
i=1\\t_i\geq
1\end{subarray}}^{l(t)}(-1)^{|t<i|+i-1}c(t-e_i,N)e^{(t)}d_A^{N-|t-e_i|-l(t)} \nonumber \\
& &\hspace{-1cm}+\sum_{s\in E_N}(-1)^{|s|+l(s)}c(s,N)e^{(s)}d_A^{N(s)+1} \label{eq5} \\
& &\hspace{-3.5cm}=\sum_{s\in E_{N+1}}\sum_{i=1}^{l(s)}\eta_i(s)(-1)^{|s<i|+i-1}c(s-e_i,N)e^{(s)}d_A^{N(s)+1}\nonumber \\
& &\hspace{-1cm}+ \sum_{s\in
E_N}(-1)^{|s|+l(s)}c(s,N)e^{(s)}d_A^{N(s)+1}. \label{eq6}
\end{eqnarray}
Putting $t=s+e_i$ in the first term of (\ref{eq4}) we obtain
(\ref{eq5}) and rewriting and changing $t$ by $s$ we obtain
(\ref{eq6}). Finally collecting similar terms in (\ref{eq3}) and
(\ref{eq6}), and using the recurrence formula we get
\[(d_A+e)^{N+1}=\sum_{s\in E_{N+1}}c(s,N+1)e^{(s)}d_A^{N(s)+1},\]
thus the proof is completed.$\blacklozenge$
\end{pf}

The following result generalizes Theorem \ref{MC2N}. It provides
an explicit formulae for the coefficients of the generalized
Maurer-Cartan equation introduced below.

\begin{thm}\label{MCMN}
We have, \[(d_{A_\R})^N=\sum_{k=0}^{N-1}c_k d_A^k,\] where
\[c_k=\sum_{\begin{subarray}{c} s\in E_N\\ N(s)=k \\ s_i<M\\ \end{subarray}}c(s,N)e^{(s)}
\hspace{.5cm}\text{and}\hspace{.5cm} c(s,N)=\sum_{\gamma\in
P_N(\emptyset,s)}\omega(\gamma).\]
\end{thm}
\begin{pf}
One checks that the coefficients $c(s,N)=\sum_{\gamma\in
P_N(\emptyset,s)} \omega(\gamma)$ satisfy the recurrence formula
of Lemma \ref{Lem23}. For this one checks that
$P_{N+1}(\emptyset,s)$ is naturally partitioned in three blocks.
The first block contains paths that are the composition of a path
$\gamma:\emptyset\to s$ in $P_N(\emptyset,s_{>1})$ with an edge
$s_{>1}\to (0,s_{>1})$ and corresponds with the first term in
(\ref{eq0re}). The second block consists of paths that are the
composition of a path $\gamma:\emptyset \to s$ in
$P_N(\emptyset,s)$ with an edge $s\to s$ and corresponds with the
second term in (\ref{eq0re}), finally the last block consists of
paths that are the composition of a path $\gamma:\emptyset \to
s-e_i$ in $P_N(\emptyset,s-e_i)$ with an edge $s-e_i\to s$ and
corresponds with the last term of (\ref{eq0re}).$\blacklozenge$
\end{pf}

Let $\Agr$ be an $M$-dga and $\Agr_\R$ an $N$-deformation over
$\R$ with $\Agr_\R=\Agr\otimes\R$. For $a\in A^1\otimes\R_+$ we
define $e_a:\Agr_\R\to\Agr_\R$ by
\[e_a(b)=ab-(-1)^{\bar{b}}ba.\]
We are assuming that the product is not graded commutative. It is
easy to see that $e_b$ is a derivation of degree 1 on
$\Agr\otimes\R_+$. Then $d_{A_\R}=d_A+e_a$ is an $N$-deformation
of $d_A$ iff $e_a$ satisfies the equation
\begin{equation}
\sum_{\begin{subarray}{c} s\in E_N \\ s_i<M\\
\end{subarray}}c(s,N)e_a^{(s)}d_A^{N-|s|-l(s)}=0. \label{eq:7}
\end{equation}
Equation (\ref{eq:7}) will be called the {\bf
$(M,N)$-Maurer-Cartan equation}. We closed this section by
formally introducing the $(M,N)$-Maurer-Cartan functor $MC_M^N(A)$
which controls deformations of the differential $d_A$ of an
$N$-dga $\Agr$.

\begin{defn}
For $N\geq M$, $a\in A^1\otimes\R_+$ is said to be an {\bf
$(M,N)$-Maurer-Cartan element} of $\Agr\otimes\R$ if $e_a$
satisfies the $(M,N)$-Maurer-Cartan equation (8). We say that $a$
is homotopic to $a'$, if $e_a$ is homotopic to $e_{a'}$ as
morphisms of $N$-dgas.
\end{defn}

\begin{defn}
We define the {\bf $(M,N)$-Maurer-Cartan} functor
$MC_M^N(A):\mathbf{Artin}\to\mathbf{Set}$ for each $M$-dga $\Agr$
over $\kr$. Functor $MC_M^N(A)$ is given by
\begin{enumerate}
\item Let $\R$ be an object of {\bf Artin}. $MC_M^N(A)(\R)$ is the set
of homotopy classes of all $(M,N)$-Maurer-Cartan elements of
$\Agr\otimes\R$.
\item If $\varphi:\R\to\R'$ be a morphism of the category
{\bf Artin} and $a$ is an $(M,N)$-Maurer-Cartan element of
$\Agr\otimes\R$, then $(1\otimes\varphi)(a)$ is an
$(M,N)$-Maurer-Cartan elements of $\Agr\otimes\R'$. Thus we obtain
a map $\varphi_{*}:MC_M^N(A)(\R)\to MC_M^N(A)(\R')$.
\end{enumerate}
\end{defn}

\smallskip

Deformation theory of $K$-dgms over an $M$-dga can be defined
similarly.

\section{Chern-Simons actions}

Let $(\Agr,m_A,d_A)$ be a 2-dga over $\kr$ and let
$(\Mgr,m_M,d_M)$ be a $2$-dgm over $(\Agr,m_A,d_A)$, consider its
$2K$-Maurer-Cartan equation, that is the equation that arises when
we deform the 2-dgm $(\Mgr,m_M,d_M)$ into a $2K$-dgm,
\(MC_{2K}(a)=(d_{End}(a)+a^2)^{K}=0,\) where $a\in End(\Mgr)$ has
degree 1. Let us assume that there exists a linear functional
$\int\!\!:End(\Mgr)\to \kr$ of degree $2K\!+\!1$, (i.e., $\int
b=0$ if $\bar{b}\neq 2K+1$) satisfying the following conditions:
\begin{enumerate}
\item $\int$ is non degenerate, that is, $\int ab=0$ for all
$a$, then $b=0$.
\item $\int d(a)=0$ for all $a$, where $d=d_{End(\Mgr)}$.
\item $\int$ is cyclic, this is $\int a_1a_2\cdots a_n=(-1)^{\bar{a_1}(\bar{a_2}\cdots\bar{a_n})}\int
a_2\cdots a_na_1$.
\end{enumerate}
We define the {\bf Chern-Simons} functional
$cs_{2,2K}:End(\Mgr)\to\kr$ by
\[cs_{2,2K}(a)=2K\!\int\!\pi (\#^{-1}(a(d_{End}(a)+a^2)^K)),\]
where
\begin{enumerate}
\item $\kr\!\!<\!\!a,d(a)\!\!>$ denotes the free $\kr$-algebra generated by
symbols $a$ and $d(a)$.
\item $\#:\kr\!<\!a,d(a)\!>\longrightarrow\kr\!<\!a,d(a)\!>$ is the linear
map defined by
\[\#(a^{i_1}d(a)^{j_1}...a^{i_k}d(a)^{j_k})=
(i_1+..+i_k+j_1+..+j_k)a^{i_1}d(a)^{j_1}.. a^{i_k}d(a)^{j_k}.\]
\item $\pi:\kr\!<\!a,d(a)\!>\longrightarrow End(\Mgr)$ is the canonical projection.
\end{enumerate}

For $K=1$ we have that $cs_{2,2}(a)$ is equal to
\[2\int\pi (\#^{-1}(a(d(a)+a^2)))=\!\! 2\int\pi (\#^{-1}(ad(a)+a^3)))
=\!\! \int ad(a)+\frac{2}{3}a^3,\] which is the Chern-Simons
functional. In general we have the following result
\begin{thm}
Let $K\geq 1$ be an integer. The Chern-Simons functional
$cs_{2,2K}$ is a Lagrangian for the $2K$-Maurer-Cartan equation,
i.e., $a\in End^1(\Mgr)$ is a critical point of $cs_{2,2K}$ if and
only if $(d(a)+a^2)^K=0$.
\end{thm}

\begin{pf}
We check that
\(\frac{\partial}{\partial\varepsilon}cs_{2,2K+2}(a+b\varepsilon)\left|_{\varepsilon=0}\right.=(2K+2)\int
bMC_{2K+2}(a).\)
\begin{eqnarray}
\frac{\partial}{\partial\varepsilon}cs_{2,2K+2}(a+b\varepsilon)\!\!&
&\!\!\left|_{\varepsilon=0}\right. \nonumber \\
=\frac{\partial}{\partial\varepsilon}(2K+2)\pi\!\!&
&\!\!\!\!\!\!\!\int
(\#^{-1}((a+b\varepsilon)MC_{2K+2}(a+b\varepsilon))\left|_{\varepsilon=0}\right. \nonumber\\
=(2K+2)\int\!\!& &\!\!\!\!\!\!\!\pi
(\#^{-1}(\frac{\partial}{\partial\varepsilon}(a+b\varepsilon)MC_{2K}(a+b\varepsilon)MC_{2}
(a+b\varepsilon)))\left|_{\varepsilon=0}\right.\nonumber\\
=(2K+2)\int\!\!& &\!\!\!\!\!\!\!\pi
(\#^{-1}(\frac{\partial}{\partial\varepsilon}(a+b\varepsilon)MC_{2K}(a+b\varepsilon))
\left|_{\varepsilon=0}\right.MC_{2}(a)\nonumber\\
+\ (2K+2)\int\!\!& &\!\!\!\!\!\!\!\pi
(\#^{-1}(aMC_{2K}(a)\frac{\partial}{\partial\varepsilon}MC_{2}(a+b\varepsilon)))\left|_{\varepsilon=0}.\right.
\label{eq8}
\end{eqnarray}
By degree reasons the second term of (\ref{eq8}) vanishes, the
inductive hypothesis yields
\begin{eqnarray*}
\frac{\partial}{\partial\varepsilon}cs_{2,2K+2}(a+b\varepsilon)\left|_{\varepsilon=0}\right.&=&
(2K+2)\int bMC_{2K}(a)MC_{2}(a)\\
&=&(2K+2)\int bMC_{2K+2}(a).\blacklozenge\\
\end{eqnarray*}
\end{pf}

For $K=2,\ 3$ the Chern-Simons functional $cs_{2,2K}(a)$ is given
by
\[cs_{2,4}(a)=\int \frac{4}{3}a(d(a))^2+2a^3d(a)+\frac{4}{5}a^5. \]
\[cs_{2,6}(a)=\int \frac{3}{2}a(d(a))^3+\frac{12}{5}a^3(d(a))^2+
\frac{6}{5}ad(a)a^2d(a)+3a^5d(a)+\frac{6}{7}a^7. \]


\subsection*{Acknowledgement}

We thank Nicol\'as Andruskiewitsch, Edmundo Castillo, Eddy
Pariguan, Sylvie Paycha and Jim Stasheff for helpful suggestions.
Thanks also to an anonymous referee for precise corrections.


\[\begin{array}{l}

\mbox{Mauricio Angel. Universidad Central de Venezuela (UCV).} \ \  \mbox{\texttt{mangel@euler.ciens.ucv.ve}} \\

\mbox{Rafael D\'\i az. Universidad Central de Venezuela (UCV).} \ \  \mbox{\texttt{rdiaz@euler.ciens.ucv.ve}} \\
\end{array}\]

\end{sloppypar}
\end{document}